\documentclass[article,12pt]{amsart}

\usepackage{a4wide,amssymb,algorithmic,algorithm,url}
\usepackage{tikz}
\usepackage[normalem]{ulem}

\usepackage{color}
\newcommand{\cg}{\color{black}} 
\newcommand{\cgb}{\color{black}} 
\newcommand{\cgr}{\color{black}}  
\newcommand{\cgg}{\color{black}} 
\newcommand{\cgp}{\color{black}}
\newcommand{\cgo}{\color{black}}

\newif\ifdetails
\detailstrue
\newcommand{\DETAIL}[1]%
{\ifdetails\par\fbox{\begin{minipage}{0.9\linewidth}\textit{Detail:}
      #1\end{minipage}}\par\fi}
\newcommand{\TODO}[1]%
{\ifdetails\par\fbox{\begin{minipage}{0.9\linewidth}\textbf{TODO:}
      #1\end{minipage}}\par\fi}

\newtheorem{lemma}{Lemma}

\newtheorem{theorem}[lemma]{Theorem}

\theoremstyle{remark}

\newcommand{\old}[1]{{}}

\title[Monte-Carlo algorithm for BPAM realization]{An algebraic Monte-Carlo algorithm for the 
Partition Adjacency Matrix realization problem}

\author{\'Eva Czabarka}
\author{L\'aszl\'o A. Sz\'ekely }
\address{\'Eva Czabarka and L\'aszl\'o A. Sz\'ekely\\ Department of Mathematics \\ 
University of South Carolina \\ Columbia, SC 29208 \\ USA}
\email{\{czabarka,szekely\}@math.sc.edu }
\thanks{The second author was supported in part by the  NSF DMS,  grant numbers 1300547 and 1600811.}
\author{Zolt\'an Toroczkai}
\address{Zolt\'an Toroczkai\\ University of Notre Dame\\
Department of Physics, Department of Computer Science and Engineering\\
225 Nieuwland Science Hall\\ Notre Dame, IN 46556\\ USA}
\email{toro@nd.edu}
\author{Shanise Walker}
\address{Shanise Walker\\ Iowa State University\\ Department of Mathematics \\
396 Carver Hall, 411 Morrill Road\\ Ames, IA 50011 \\USA}
\email{Shanise1@iastate.edu}
\subjclass[2010]{Primary 05D15; secondary 05D40; 05C82; 68W20}
\keywords{perfect matching, exact matching, degree sequence, bipartite graph, Hall's Theorem,  
Joint Degree Matrix, Partition Adjacency Matrix, Monte-Carlo algorithm, computational complexity}

\begin{document}

\begin{abstract} 
{\cgb The graphical realization of a given degree sequence and given partition} adjacency 
matrix simultaneously is a relevant problem {\cgb in data driven modeling of networks}. Here
we formulate common generalizations of this problem and the Exact Matching 
Problem, and solve them with an algebraic Monte-Carlo algorithm {\cgb that} runs in polynomial 
time if the number of partition classes is bounded. 
\end{abstract}

\maketitle

\section{Introduction}
{\cgb In data driven modeling of complex networks one often needs 
to sample from ensembles of graphs that 
share characteristics with an observed network. These characteristics act as constraints for the sampling
procedure and they may be reproduced exactly (``sharp constraints'') in every sampled graph 
or in expected value (average constraints) over the ensemble. The most natural such characteristic is the 
degree sequence. The degree sequence, however, has many graphical realizations in general, with varying 
properties, e.g.,} either showing assortativity or disassortativity (the extent to which vertices of similar 
degrees are connected or not). For example, social networks tend to be assortative, while biological and 
technical networks tend to be disassortative. {\cgb Thus, in order to model such situations, one also has to specify the degree correlations. 
The simplest way of achieving that is via providing the so-called Joint Degree Matrix (JDM), whose entries are the number of
edges between degree $i$ and degree $j$ vertices, for all $i$ and $j$ degree values.} Note that the JDM also
specifies the degree sequence itself, uniquely \cite{dutle}. The JDM {\cgo received} considerable attention in the literature \cite{aman1,aman2,dutle,jdmMCMC,sadeghis,patrinos,rinaldo,sadeghirinaldo,stantonpinar} and it is well 
understood \cite{aman1,aman2,dutle,patrinos,stantonpinar}. {\cgb Reference \cite{JDMsampling} presents an 
exact algorithm for constructing simple graphs with a prescribed JDM.} 

{\cgo  However, to model real world networks
even JDM level constraints are not always sufficient. In particular,  Orsini et al. \cite{Orsini} demonstrate this on several networks including the
Internet (autonomous systems level), the distributed PGP cryptosystem, US airport network, protein interactions, brain fmri functional networks and an English word adjacency network. To capture most of the ``signal" found in the structure of a real-world network they also had to include correlations beyond degree-degree correlations, such as clustering coefficients, i.e., small subgraph correlations, called collectively as $dk$-series (1$k$ is purely degree distribution, 2$k$ is joint degree distribution, etc). When trying to generate and sample graphs with these prescriped correlations, these authors observe that already at the $d=3$ level the process does not converge and modelling essentially fails. As briefly discussed in \cite{Orsini}, this graph generation process can be described as a Boolean constraint 
satisfaction SAT problem, in which the variables are the elements of the adjacency matrix whose values need to be set (to 0 or 1) such that a set of constraints expressed in terms of functions of the marginals (degrees) are satisfied. From this point of view it is thus expected that the problem eventually becomes NP-complete (3-SAT), which indeed was experienced in \cite{Orsini} through the failure of the algorithms to converge.} 

{\cgo The graph construction problems above all have their constraints related to some structural properties of the graph. However, in many real world situations
there are also externally imposed constraints, such as group membership that is not modelled by the approaches above. For example, one might study a 
network at different levels of resolution: we may look at a large organization as a network of interactions between teams or departments but also at the connections between the individuals throughout the organization and ask questions related to the performance of the organization as a whole as function of these networks. One can certainly think of the teams and departments as a {\em partitioning} of the individuals into groups/classes and the connections between them as a partitioning of the edges. }{\cgb In 2014 the first author introduced the concept of the Partition Adjacency Matrix (PAM), in order to also accomodate such, more general classes of constraints in network modeling {\cgo \cite{GRAPHS}}.} Given a simple graph and a partition
of its vertex set, entries of the {\cgb PAM} count {\cgb the} number of edges {\cgb between, and within} the partition classes. If the partition consists of singleton sets of vertices, then it specializes to the familiar Adjacency Matrix of the graph, while the {\cgb JDM} of a graph is {\cg a special} 
PAM, with all classes composed of vertices having the same degree.

{\cgb In contrast with the JDM, however, much less is known about PAMs, which are, as explained above, an important notion in data driven modeling of networks.}
{\cgg While a JDM determines the degree sequence of a realizing graph, a PAM does not.}
{\cgb
Similarly to JDM problems, PAM problems include existence
({\em Is there a simple graph with a {\cgg given degree sequence and} given PAM?}), construction ({\em Provide an algorithm that constructs such simple graphs\cgr!}), sampling ({\em Provide an algorithm that can sample such graphs with prescribed distribution\cgr!}) and counting
problems ({\em How many simple graphs are there to realize a {\cgg given degree sequence and } given PAM?}), in increasing order of their difficulty. 
Here we focus on the realization and construction problems.} 

{\bf Partition Adjacency Matrix realization problem}:
Given a  set $W$ and natural numbers $d(w)$ associated with $w\in W$,  a $W_i:i\in I$ partition of $W$, and
natural numbers $c(W_i,W_j)$
associated with unordered pairs of partition classes, is there  a
simple graph on the vertex set $W$ with degree  $d(w)$  for every $w\in W,$ and with exactly $c(W_i,W_j)$ edges with endpoints in
$W_i$ and $W_j$?

{\bf  Partition Adjacency Matrix construction problem}: {\cgp Construct} such a graph, if the answer to the realization problem is {\cgb affirmative}. 

{\cgp Reference} \cite{skeleton} conjectures that the realization problem is NP-complete, 
and {\cgb here we also support} this conjecture. 
The {\em skeleton} of a PAM is  
the graph, whose
vertices are the partition classes, and two partition classes, $W_i$ and $W_j$ are joined by an edge, if $c(W_i,W_j)>0$. {\cgp Reference} \cite{skeleton}
found polynomially solvable instances of the realization problem for {\cgb two partition classes ($|I|=2$)}, and also for loopless unicyclic skeleton graphs. 
In the Bipartite PAM problem the skeleton graph is bipartite and loopless. 

A stronger version {\cgb of the problems above is when there is also a forbidden subgraph that all graphical realizations must avoid. Such problems arise {\cgg in part} for algorithmic reasons in direct construction algorithms that add edges sequentially: the existing edges forbid the addition of further edges between the same pairs of vertices 
in the graph being constructed \cite{kim}, \cite{JDMsampling}. Thus, we formulate:}

 {\bf Partition Adjacency Matrix realization/construction problems in the presence of a blue graph}: {\cgp In addition to the contraints of the PAM realization problem,} 
 a graph $B$  (the blue graph) is given on the vertex set $W$. Is there a realization {\cgb that is} not using any edges from $B$? If yes, construct such a graph.
 
To simplify {\cgb the} discussion, we assume that for any PAM realization/construction problem, the obvious and easy-to-check necessary conditions that $d$ is the degree sequence {\cgb of a simple graph}, and {\cgb that}
$$\sum_{v\in W_i} d(v) =c(W_i,W_i)+\sum_j c(W_i,W_j),$$
$$\sum_v d(v) = \sum_i c(W_i,W_i) +\sum_i \sum_j c(W_i,W_j)$$
hold. 
In an earlier version of this manuscript \cite{JDM-MC}, we gave an algebraic Monte-Carlo algorithm for the blue graph version of  
Bipartite PAM realization  problem. This algorithm runs
in polynomial time if the number of partition classes is bounded. We are indebted to Andr\'as Frank (Budapest), who kindly called our attention to the analogous Exact Matching Problem \cite{website}, and asked if the two problems admit a common generalization, {\cgr i.e. a third problem, of which the first two problems are specific instances, such that third problem allows algorithmic solution like the first two.} {\cgb Let us recall  the}

{\bf Exact Matching Problem}: Given a  graph $G$, whose edges {\cgb are} colored red or green, is there a  perfect matching with exactly $m$ red edges {\cgo in the matching}?

The Exact Matching Problem originates from Papadimitriu and Yannakakis \cite{exactprobl}, and Lov\'asz  proposed a Monte-Carlo algorithm for it. {\cgr Lov\'asz' algorithm, which he never published, is based on the general ideas in his paper \cite{lovasz}, and is described 
{\cgo by Mulmuley, Vazirani and Vazirani} in \cite{vazirani} pp. 111.}
No deterministic polynomial time algorithm is known for the  Exact Matching Problem.
{\cgb Here we provide} {\cgr the promised} common generalization:

{\bf Dominating f-factor Problem}:
Given a graph $G$ on $n$ vertices, disjoint subsets $E_1,...,E_k\subset E(G)$, integers $m_1,...,m_k$, and prescribed degrees $d_1,...,d_n$ associated
with the  vertices $v_1,...,v_n$ of $G$, is there a subgraph $G'$ of $G$, such that $v_i$ has degree $d_i$ in $G'$ for all vertices,
and $G'$ has at least $m_j$ edges from the edge set $E_j$, for {\cgo all} $j=1,...,k$?

{\bf Dominating Matching Problem}: Given a graph $G$, disjoint subsets $E'_1,...,E'_k\subset E(G)$, integers $m_1,...,m_k$, is there  a
perfect matching in $G$, which uses at least   $m_j$ edges from the edge set $E'_j$, for {\cgo all} $j=1,...,k$?

Clearly the Dominating Matching Problem is a special case of the Dominating f-factor Problem, {\cgo where} every degree is one.
In Section~\ref{sec:tutte}, we will show  using the Tutte gadget that   the Dominating f-factor Problem can be solved through solving 
a Dominating Matching Problem on about $n^2$ vertices. The Exact Matching problem is an instance of the Dominating Matching Problem,
where $k=2$, $E_1$ is the set of red edges, $E_2$ is the set of green edges, $m_1=m$, $m_2=|E(G)|-m$. The PAM realization problem
is an instance of the Dominating f-factor Problem in the following way: $G$ is the complement of the blue graph $B$, the disjoint edge subsets
are $E_{ij}=\{\{u,v\}\in E(G):u\in W_i, v\in W_j\}$ for $ i\leq j$ (assuming without loss of generality that $I$ is an ordered set), $m_{ij}=c(W_i,W_j)$ for $i<j$ and 
$m_{ii}=c(W_i,W_i) $.

{\cgb Here} we provide an algebraic Monte-Carlo algorithm for the  Dominating Matching Problem, and hence for the  Dominating f-factor Problem, which runs
in polynomial time under the assumption that 
$$\prod_{i\in I} (m_{i}+ 1)=   O(polynomial(n)).$$ 
This assumption certainly holds
if $|I|$  {\cgb stays} bounded, while $n$ {\cgb grows}. If the algorithm returns  {\tt TRUE}, then the sought-after graph exists, if the algorithm
returns {\tt FALSE}, then with high probability (whp) such a graph does not exist. The correctness of the algorithm hinges on the 
Schwartz-Zippel Lemma~\cite{Schwartz}.  The realization 
algorithm and its correctness {\cg are} described in Section~\ref{sec:main}. We will also conclude that constructing an actual solution
is not harder than the decision problem. We conclude the paper with some complexity results in Section~\ref{concluding}.


\section{The Tutte gadget}\label{sec:tutte}

{\cgb Clearly, the}  standard degree sequence realization problem is a relaxation of the PAM  problem, where 
 we do not care for satisfying the $c_{ij}$ conditions. Havel \cite{havel} and Hakimi \cite{hakimi} solved the degree sequence realization problem and
 Ryser \cite{ryser} solved the bipartite degree sequence realization problem. {\cgb We next}
 use a result of Tutte \cite{gadget} to connect degree sequence realization to the existence of {\cgb a} perfect matching in a bigger graph, the {\it Tutte gadget}. 

 Initially, we are given a degree sequence realization problem on the vertices in $V$, i.e. for each $v\in V$ we are given a proposed degree $d(v)$. We are also given a set of blue -- or forbidden -- edges $B$ that our realization is not allowed to use. For a vertex $v\in V$, {\cgb let} $N_B(v)=\{u:\{u,v\}\in B\}$
{\cgb denote} the set of blue neighbors of $v$, and
 $S_v= V\setminus\left(\{v\}\cup N_B(v)\right)$ {\cgb denote} the set of allowed neighbors. Without  restrictions in the degree
 sequence realization problem, we have $S_v=V\setminus\{v\}$. 
 The setup of this problem implies that $u\in S_v$ iff $v\in S_u$, and we will also assume further that 
 for each $v\in V$ $|S_v|\ge d(v)$ holds (otherwise a realization obviously {\cg cannot} exist).

The Tutte gadget of the degree sequence realization problem with a set $B$ 
of blue edges is a graph $T$ such that 
$$V(T)=\{v^u:v\in V,u\in S_v\}\cup \{a^v_1,...,a^v_{|S_v|-d(v)}:v\in V   \}$$ and 
$$
 E(T)=\left\{\{v^u,u^v\}: v\in V, u\in S_v\right\}\cup\left\{\{v^u, a^v_i\}: v\in V, u\in S_v, i=1,2,...,|S_v|-d(v)\right\}.
$$

For  the degree sequence realization problem, with $B=\emptyset$, for each $v\in V$, $S_v=V\setminus\{v\}$ and the degree condition becomes
$d(v)\le |V|-1$. The Tutte gadget is a graph with $2n(n-1)-\sum_v d(v)$ vertices.

\old{
In the case of the bipartite degree sequence problem, when we have vertex classes $U,W$ with $|U|=m$ and $|W|=n$, and edges
only allowed between $U$ and $W$,
the Tutte gadget becomes
$$V(T)=\{w^u,u^w:w\in W,u\in U\}\cup
 \{a^w_1,\ldots,a^w_{m-d(w)}:w\in W \}\cup\{a^u_1,\ldots,a^u_{n-D(u)}:u\in U\}$$ and 
\begin{eqnarray*}
 E(T)&=&\left\{\{w^u,u^w\}: w\in W, u\in U\right\}\cup\left(\bigcup_{w\in W,u\in U}\left\{\{w^u, a^w_i\}:i\in[m-d(w)]\right\}\right)\cup\\
 & &\,\,\,\left(\bigcup_{u\in U,w\in W}\left\{\{u^w, a^u_i\}:i\in[n-D(u)]\right\}\right),\end{eqnarray*}
where for each $u\in U${\cg,} $D(u)\le n$ and for each $w\in W$, $d(w)\le m$. 
For the bipartite degree sequence problem the Tutte gadget $T$ is a bipartite graph with partite classes of size $2nm-{E}$. 
}



The Tutte gadget is relevant for the following property:
{\cgb it} has a perfect matching if and only if a  graph solves the corresponding degree sequence realization problem; 
 furthermore, if some $\{w^u,u^w\}$ edges are present in the perfect matching, then  the corresponding $\{w,u\}$ edges provide a graph solving this degree sequence realization problem, and if some $\{w,u\}$ edges provide a 
 graph solving the  degree sequence realization problem, then the corresponding $\{w^u,u^w\}$ edges in $T$ are part of a perfect matching of $T$. This property is well-known and is also easy to verify.
 
 Furthermore, if an edge $\{u,v\}\in E(G)$ belongs to an edge set $E_i$, put the edge $\{w^u,u^w\} \in E(T)$ into $E'_i$, when solve the 
 Dominating f-factor Problem from the Dominating Matching Problem using the Tutte gadget.

\section{The Dominating Matching Problem}\label{sec:main}

Let $A$ be a skew-symmetric matrix, i.e. $A=-A^T$, and assume that $A$ has an even order $2n$. 
The Pfaffian of $A$ is defined as
$${\rm Pf}(A)=\sum_\pi {\rm sign}(\pi)\cdot {}_{i_1}[A]_{j_1}\cdot {}_{i_2}[A]_{j_2}\cdots{}_{i_n}[A]_{j_n},$$
where  $\pi$ runs through permutations of the form
$$\pi=\begin{pmatrix}1& 2 & 3 & 4& ...& 2n-1 &2n \\ i_1& j_1 &i_2 & j_2 &... & i_n & j_n \end{pmatrix}$$
under the assumptions $i_1<j_1, i_2<j_2,\ldots i_n<j_n$ and $i_1<i_2<\cdots <i_n$, and ${\rm sign}(\pi)=\pm 1$, the sign of the
permutation $\pi$. {\cgr For more background on the Pfaffian, see \cite{lovaszplummer}.} 
Cayley \cite{cayley} {\cgr and Muir \cite{Muir1,Muir2}} proved that $({\rm Pf}(A))^2=\det(A)$. Note that the summation for $\pi$ can be thought of as a 
summation over the perfect matchings of $2n$ elements.

Assume  now that we are given a graph $G$ for the Dominating Matching Problem. We will assume that the graph has an even number of vertices, say $2n$, otherwise it cannot have a perfect matching. Fix an arbitrary orientation $\vec G$ of the graph $G$. For the arc 
$i\rightarrow j$ introduce a variable $x_{ij}$, and define $A$ by
$${}_i[A]_j=\begin{cases} x_{ij}  & \text{if\ }   i\rightarrow j \\ -x_{ij}  & \text{if\ }  j\rightarrow i \\ 0 & \text{otherwise}. \end{cases}$$
The variables $x_{ij}$ are independent of each other. It is clear that $G$ has a perfect matching if and only if the {\em polynomial}
${\rm Pf}(A)$ is not termless, i.e not the zero polynomial, as {\cgb cancellation} of terms is not possible. Tutte's theorem \cite{skew}, that $G$ has a perfect
matching if and only if the polynomial $\det (A)$ is not the zero polynomial follows from Cayley's theorem. Introduce now 
additional new variables, $z_\ell$ associated with the edge set $E'_\ell$, for $\ell=1,2,...,k$. Define the matrix $A^*$ by the substitutions
$x_{ij}\leftarrow x_{ij}z_\ell$ for all $\{i,j\}\in E'_\ell$ in $A$, and not changing $x_{ij}$ if $\{i,j\}\notin \cup_{\ell=1}^k E'_\ell$. A matching 
that defines a term in ${\rm Pf}(A)$ solves the Dominating Matching Problem if and only if for every $\ell$, the exponent of $z_\ell$ is
at least $m_\ell$, for $\ell=1,2,...,k$.


Now we need some properties of the  difference operator acting on multivariate polynomials. For a polynomial $f(x,y,z,...)$, set 
$$\triangledown_x f= f(x,y,z,...)- f(x-1,y,z,...).$$
We will use {\cg products} of these operators to indicate juxtaposition, and consequently 
 $\triangledown_x^k$ will denote the repetition of the operator $\triangledown_x$ $k$ times. $\triangledown_x^0$ is the identity operator. Note that unless $f$ is identically zero, applying
$\triangledown_x$ strictly decreases the degree of $x$ in the polynomial. Therefore, if the degree of $x$ in $f$ is less than $k$, then $\triangledown_x^k f$ is identically 0, and  $\triangledown_x^k x^k=k!   \not= 0$.
It is well-known that
$$\triangledown_x^k f(x)=\sum_{\ell=0}^k (-1)^\ell \binom{k}{\ell}f(x-\ell).
$$
Furthermore, as 
$$\triangledown_x \triangledown_y f=f(x,y,z,...)- f(x-1,y,z,...)-f(x,y-1,z,...)+f(x-1,y-1,z,...) =  \triangledown_y \triangledown_x f ,        $$
the order of $\triangledown $ operators associated with different variables is freely interchangeable. 
For any function $f$ in variables $z_{\ell}$, and possibly other variables not shown, the following iterated difference, which is put into product notation,  can be computed formally:
\begin{equation*} \Biggl( \prod_{\ell=1}^k    \triangledown_{z_{\ell}}^{m_{\ell} }\Biggl)  f(z_1,\ldots, z_{\ell}, \ldots, z_k)= 
\end{equation*}
\begin{equation}
\sum_{u_{1}=0}^{m_{1}}\cdots\sum_{u_{\ell}=0}^{m_{\ell}}
\cdots\sum_{u_{k}=0}^{m_{k}}
\left(\prod_{\ell=1}^k(-1)^{u_{\ell}}\binom{m_\ell}{u_\ell}\right)f(z_{1}-u_{1},\ldots,  z_{\ell}-u_\ell   ,\ldots z_{k}-u_k). \label{ugly}
\end{equation}

We are ready to claim the key fact behind our algorithm: the polynomial
\begin{equation*}
 \Biggl( \prod_{\ell=1}^k    \triangledown_{z_{\ell}}^{m_{\ell} }\Biggl)  {\rm Pf}\left(A^{\star}(z_1,\ldots,z_k)\right)=
\end{equation*}
{\cgg\begin{equation}\label{detB3}
\sum_{u_{1}=0}^{m_{1}}\cdots\sum_{u_{\ell}=0}^{m_{\ell}}
\cdots\sum_{u_{k}=0}^{m_{k}}
\left(\prod_{\ell=1}^k(-1)^{u_{\ell}}\binom{m_\ell}{u_\ell}\right){\rm Pf}\left(A^{\star}(z_{1}-u_{1},\ldots,  z_{\ell}-u_\ell   ,\ldots z_{k}-u_k)\right)
\end{equation}}
is not identically 0 if and only if the Dominating Matching Problem has a solution, as 
no monomial can be  a multiple of another. 
{\cg Thus,} the Dominating Matching Problem
 boils down to checking whether the polynomial (\ref{detB3}) is identically zero or not. Make random substitutions into all 
variables of (\ref{detB3}), if this polynomial is not identically zero, then whp after a number of substitutions we obtain a nonzero value. In this 
case the answer to the  problem is a (correct) {\tt TRUE}. If we always get zero values, the answer returned is a {\tt FALSE}, and it is correct 
whp. (We give a more detailed analysis below.) From a computational point of view, the issue is whether we can compute substituted values
of (\ref{detB3}) in polynomial time. Note that the Pfaffian with integer entries (or with entries from an integral domain) can be evaluated
efficiently, similarly to the evaluation of a determinant \cite{efficient}.

While the polynomial ${\rm Pf}(A^{\star})$ is not computable in polynomial time, the result of substituting numbers into all variables {\cgb is}. 
Indeed, (\ref{detB3}) expanded in (\ref{ugly}) with     $ f(z_1,\ldots, z_{\ell}, \ldots,z_k)={\rm Pf}(A^{\star})                           $ is just a weigthed sum of values of ${\rm Pf}(A^{\star})$ after the substitutions
$ z_{\ell}\leftarrow  z_{\ell}-u_{\ell}$ $ (\ell=1,2,...,k)$ for every $0\leq u_\ell \leq m_\ell$. In other words, for every attempt to 
substitute random numbers, we have to evaluate $\prod_{\ell=1}^k (m_\ell+ 1)$ numerical Pfaffians, a polynomial number of steps in $n$.

Recall the Schwartz-Zippel Lemma \cite{Schwartz}, where non-zero polynomial means that at least one term comes with nonzero coefficient.
\begin{lemma} For a field ${\mathbb F}$,
let $f\in {\mathbb F}[x_1,x_2,...,x_t] $ be a non-zero polynomial of degree $d$ and $\Omega \subseteq  {\mathbb F}$ a finite set,
$|\Omega |=N$. Let $Z(f,\Omega)$ denote the set of roots from $\Omega^n$, i.e. 
$$Z(f,\Omega)= \{(\alpha_1,\alpha_2,...,\alpha_t)\in\Omega^n: f(\alpha_1,\alpha_2,...,\alpha_t)=0\}.$$
Then $|Z(f,\Omega)|\leq dN^{t-1}$, and the probability that $f$ vanishes on randomly and independently selected uniformly random elements 
of $\Omega$ is at most $d/N$.
\end{lemma}
The polynomial (\ref{detB3}) has degree at most $2n$. Let $p$ be a prime number, such that
$p\geq    2n^2$. One can find such a prime using Bertrand's Postulate (better estimates on gaps between primes exist) and prime testing 
the numbers one after the other. Set ${\mathbb F}=\Omega=GF(p)$. We compute (\ref{detB3}) in $GF(p)$, i.e. we do the calculations
mod $p$. Note that the polynomial (\ref{detB3}) is non-zero over $GF(p)$ as well if a solution to the Dominating Matching Problem exists, since after {\cgb taking the derivatives} we get coefficients at the terms
that are products of numbers at most $2n$.

Substituting randomly and uniformly selected  elements of  $GF(p)$ into the variables of $A^{\star}$ and its {\cg translates}, the probability of getting
a 0 value for the expression (\ref{detB3}) if it is not the identically 0 polynomial, is at most $2n/N=2n/p\leq 1/n$, according to the Lemma. 


\begin{theorem}\label{state}
There is a Monte-Carlo algorithm for the Dominating Matching Problem and the Dominating f-factor Problem, which runs
in polynomial time under the assumption that 
$\prod_{i\in I} (m_{i}+ 1)=   O(polynomial(n)),$ 
which  certainly holds
if $|I|$  stays bounded. If the algorithm returns  {\tt TRUE}, then the sought {\cg after} graph exists, if the algorithm
returns {\tt FALSE}, then with high probability (whp) such a graph does not exist.
\end{theorem}

An actual solution easily can be found by testing iteratively whether an edge can be included in the matching in a modified problem,
a standard approach \cite{vazirani,kim}. In the first version of this manuscript \cite{JDM-MC} we provided {\cgb a} pseudocode for the Bipartite PAM realization and construction problems.

\begin{theorem}\label{state1}
There is a Monte-Carlo algorithm to construct a solution for the Dominating Matching Problem or the Dominating f-factor Problem, which runs
in polynomial time under the assumption that 
$\prod_{i\in I} (m_{i}+ 1)=   O(polynomial(n)),$ 
which  certainly holds
if $|I|$  stays bounded. If the algorithm returns  a construction, then it is {\cgo a} correct solution, and if a correct solution exists, a construction
is found whp.
\end{theorem}

\section{Concluding remarks}\label{concluding}

Our algorithm for the Dominating Matching Problem, if specialized for the Exact Matching Problem, is  different from
from Lov\'asz' algorithm \cite{vazirani}. 
We believe, {\cgb however,} that the same techniques {\cgb may} {\cgr also} be used  to solve the 
Dominating Matching Problem.

{\cg Here we did not} attempt to optimize and estimate the running time of the algorithms, as they are very far from practical. 
We repeat here that no deterministic polynomial time algorithm is known for the Exact Matching problem, not even for bipartite graphs.
Hence no deterministic polynomial time algorithm is known for the Dominating f-factor and Dominating Matching Problems.

We are thankful to Stefan Lendl (Graz) for {\cgb bringing} to our attention reference \cite{plaisted}. Ref. \cite{plaisted} shows that given a bipartite 
graph $G$ and  a partition  $V_1,...,V_s$ and $U_1,...,U_\ell$   of the paritite classes, the decision problem whether a perfect matching $M$ exists with at most 1 edge between any
pair of partition classes is NP-complete. It is easy to see that this problem is equivalent to the following instance of the Dominating
f-factor Problem: the graph is $G$, the prescribed degree is $d_G(v)-1$ for vertex $v$, the $E_{ij}$ edge sets are $E(G)\cap (V_i\times U_j)$,
and $m_{ij}= |E(G)\cap (V_i\times U_j)|-1$. Hence the  Dominating
f-factor Problem is also NP-complete. 

Ref. \cite{geerdes} noted that 3-dimensional  perfect matching problem  in 3-partite graphs can be reduced to the problem of 
finding a multicolored perfect matching in an $n$-colored bipartite graph $K_{n,n}$. This gives another proof for the fact that the
Dominating Matching Problem is NP-complete.
{\cgo Reference} \cite{skeleton} conjectures that the PAM realization problem (with empty blue graph) is already NP-complete.


\end{document}